\documentstyle[amssymb]{amsart}
\newcommand{\de}{\delta}
\newcommand{\al}{\alpha}
\newcommand{\si}{\sigma}
\newcommand{\tha}{\theta}
\newcommand{\pa}{\partial}

\begin{document}
\begin{titlepage}
\begin{flushright}
M\'exico ICN-UNAM 97-14
\end{flushright}
\vskip 8mm
\begin{center}
{\LARGE  Lie Algebras in Fock Space}

\vskip 0.3cm
{\large by}
\vskip 0.3cm
{\large A. Turbiner}$^{\star}$
\vskip .2cm
{\it Instituto de Ciencias Nucleares, UNAM, \\Apartado Postal 70-543,
04510 Mexico D.F., Mexico}
\vskip 0.3cm
\end{center}

\begin{quote}
A catalogue of explicit realizations of representations of Lie (super)
algebras and quantum algebras in Fock space is presented.
\end{quote}

\vskip 5cm

\noindent
$^{\star}$On leave of absence from the Institute for Theoretical
and Experimental Physics,
Moscow 117259, Russia\\E-mail: turbiner@@axcrnb.cern.ch,
turbiner@@roxanne.nuclecu.unam.mx

\end{titlepage}

This article is an attempt to present a catalogue of known representations
of (super) Lie algebras and quantum algebras acting on different Fock spaces.
Of course, we do not have as ambitious goal as to present a complete list of all
representations of all possible algebras, but we plan to present some
of them, those we consider important for applications, mainly restricting
ourselves by those possessing finite-dimensional representations.  Many
representations are known in the folklore spread throughout the literature
under offen different names\footnote{For instance, in nuclear physics some
of them are known as boson representations}. Therefore, we provide
references according to our taste, knowledge and often quite arbitrarily.
This work does not pretend to be totally original.
Throughout the text we usually consider complex algebras.

\bigskip
\begin{center}
{\Large Lie Algebras}

{\it 1.} $sl_2$-algebra.
\end{center}

\renewcommand{\theequation}{A.1.{\arabic{equation}}}
\setcounter{equation}{0}

Take two operators $a$ and $b$ obeying the commutation relation
\begin{equation}
\label{e1.1}
             [a,b] \equiv ab  -  ba \ =\ 1,
\end{equation}
with the identity operator on the r.h.s. -- they span
the three-dimensional Heisenberg algebra. By definition the universal enveloping
algebra of the Heisenberg algebra is the algebra of all normal-ordered
polynomials in $a,b$: any monomial is taken to be of the form $b^k a^m$
\footnote{Sometimes this is called the Heisenberg-Weil algebra}. If, besides
the polynomials, we also consider holomorphic functions in $a,b$,
the extended universal enveloping algebra of the Heisenberg
algebra appears. The Heisenberg-Weil algebra possesses the internal automorphism,
which is treated as a certain type of quantum canonical transformations
\footnote{This means that
there exists a family of the elements of the Heisenberg-Weil algebra obeying the
commutation relation (\ref{e1.1})}.
We say that the (extended) Fock space appears if we take the (extended)
universal enveloping algebra of the Heisenberg algebra and add to it
the vacuum state $|0>$ such that
\begin{equation}
\label{e1.2}
a|0>\  = \ 0\ .
\end{equation}
One of the most important realizations of (\ref{e1.1}) is the
coordinate-momentum representation:
\begin{equation}
\label{e1.3}
a\ =\ \frac{d}{dx} \equiv \pa_x\ ,\ b\ =\ x\ ,
\end{equation}
where $x \in {\bf C}$. In this case the vacuum is a constant, say,
$|0>\  = \ 1$. Recently a finite-difference
analogue of (\ref{e1.3}) has been found \cite{st},
\begin{equation}
\label{e1.4}
a\ =\ {\cal D}_+ ,\ b\ =\ x(1-\de{\cal D}_-)\ ,
\end{equation}
where
\[
{\cal D}_+ f(x) = \frac{f(x+\de) - f(x)}{\de}\ ,
\]
is the finite-difference operator, $\de \in C$ and
${\cal D}_+\rightarrow {\cal D}_-$, if $\de \rightarrow -\de$.

{\bf (a).}
It is easy to check that if the operators $a,b$ obey (\ref{e1.1}), then the
following three operators
\[
J^+_n = b^2 a - n b\ ,
\]
\begin{equation}
\label{e.a.1}
J^0_n = ba - {n \over 2}\ ,
\end{equation}
\[
J^-_n=a\ ,
\]
span the $sl_2$-algebra with the commutation relations:
\[
[J^0,J^{\pm}]=\pm J^{\pm}\ ,\  [J^+,J^-]=-2J^0\ ,
\]
where $n \in C$.
For the representation (\ref{e.a.1}) the quadratic Casimir operator is equal to
\begin{equation}
\label{e.a.2}
C_2 \equiv \frac{1}{2}\{J^+,J^-\} - J^0 J^0 = -\frac{n}{2}
\bigg(\frac{n}{2}+\frac{1}{2}\bigg)\ ,
\end{equation}
where $\{\ ,\ \}$ denotes the anticommutator.
If $n$ is a non-negative integer, then (\ref{e.a.1}) possesses
a finite-dimensional, irreducible representation  in the Fock space leaving
invariant the space
\begin{equation}
\label{e.a.3}
{\cal P}_{n}(b) \ = \ \langle 1, b, b^2, \dots , b^n \rangle ,
\end{equation}
of dimension $\dim{\cal P}_{n}=(n+1)$.

Substitution of (\ref{e1.3}) into (\ref{e.a.1}) leads to a well-known
representation of $sl_2$-algebra of differential operators of the first order
\footnote{This representation was known to Sophus Lie.}
\[
J^+_n = x^2 \pa_x - n x \ ,
\]
\begin{equation}
\label{e.a.4}
J^0_n = x \pa_x - {n \over 2}\ ,
\end{equation}
\[
J^-=\pa_x\ ,
\]
where the finite-dimensional representation space (\ref{e.a.3})
becomes the space of polynomials of degree not higher than $n$
\begin{equation}
\label{e.a.5}
{\cal P}_{n}(x) \ = \ \langle 1, x, x^2, \dots , x^n \rangle \ .
\end{equation}

{\bf (b).}
The existence of the internal automorphism of the extended universal enveloping
algebra of the Heisenberg algebra, i.e $[\hat a (a,b),\hat b(a,b)]=[a,b]=1$
allows to construct different representations of the algebra $sl_2$
by $a\rightarrow \hat a, b \rightarrow \hat b$ in (\ref{e.a.1}).
In particular, the internal automorphism of the extended universal enveloping
algebra of the Heisenberg algebra is realized by
the following two operators,
\[
\hat a=\frac{(e^{\de a} -1)}{\de}\ ,
\]
\begin{equation}
\label{e.b.1}
\hat b= b e^{ - \de a}\ ,
\end{equation}
where $\de$ is any complex number. If $\de$ goes to zero
then $\hat a \rightarrow a, \hat b \rightarrow b$.
In other words, (\ref{e.b.1}) is a 1-parameter quantum canonical
transformation of the deformation type of the Heisenberg algebra (\ref{e1.1}).
It is the quantum analogue of a point canonical transformation. The substitution
of the representation (\ref{e.b.1}) into (\ref{e.a.1}) results in the
following representation of the $sl_2$-algebra
\[
J_n^+= ({b \over \de}-1)b e^{ - \de a}(1 - n - e^{ - \de a})\ ,
\]
\begin{equation}
\label{e.b.2}
J_n^0= \frac{b}{\de} (1 - e^{ - \de a}) -\frac{n}{2}\ ,
\ J^-= \frac{1}{\de} ( e^{  \de a}-1)\ .
\end{equation}
If $n$ is a non-negative integer, then (\ref{e.b.2}) possesses
a finite-dimensional irreducible representation of the dimension
$\dim{\cal P}_{n}=(n+1)$ coinciding with (\ref{e.a.3}). It is worth noting
that the vacuum for (\ref{e.b.1}) remains the same, for instance (\ref{e1.2}).
Also the value of the quadratic Casimir operator for (\ref{e.b.2})
coincides with that given by (\ref{e.a.2}).

The operator $\hat a$ in the particular representation (\ref{e1.4})
becomes the well-known translationally-covariant finite-difference operator
\begin{equation}
\label{e.b.3}
\hat a f(x) = \frac{(e^{\de \pa_x} -1)}{\de} f(x)\ =\ {\cal D}_+ f(x)
\end{equation}
while  $\hat b$ takes the form
\begin{equation}
\label{e.b.4}
\hat b f(x) = x e^{-\de \pa_x} f(x)\ =
\ x f(x-\de) = x(1-\de{\cal D}_-) f(x) \ .
\end{equation}
After substitution of (\ref{e.b.3})--(\ref{e.b.4}) into (\ref{e.b.2}) we arrive at
a representation of the $sl_2$-algebra by finite-difference operators,
\[
J^+_n= x({x \over \de}-1) e^{ - \de \pa_x} (1 - n -
e^{ - \de \pa_x})\ ,
\]
\begin{equation}
\label{e.b.5}
J^0_n= {x \over \de} (1 - e^{ - \de \pa_x})-\frac{n}{2}\ ,
\ J^-= {1 \over \de} ( e^{\de \pa_x}-1)\ ,
\end{equation}
or, equivalently,
\[
J^+_n= x(1-\frac{x}{\de})(\de^2 {\cal D}_-{\cal D}_- -
(n+1)\de{\cal D}_- + n)\ ,
\]
\begin{equation}
\label{e.b.6}
J^0_n= x {\cal D}_- - \frac{n}{2}\ ,
\ J^-= \ {\cal D}_+ .
\end{equation}
The finite-dimensional representation space for (\ref{e.b.5})--(\ref{e.b.6})
for integer values of $n$ is again given by the space (\ref{e.a.5})
of polynomials of degree not higher than $n$.

{\bf (c).}
Another example of quantum canonical transformation is given by the
oscillatory representation
\[
\hat a = \frac{b+a}{\sqrt{2}}\ ,
\]
\begin{equation}
\label{e.c.1}
\hat b = \frac{b-a}{\sqrt{2}}\ .
\end{equation}
Inserting (\ref{e.c.1}) into (\ref{e.a.1}) it is easy to check that
the following three generators form a representation of the
$sl_2$-algebra,
\[
J^+_n = \frac{1}{2^{3/2}}[b^3 + a^3 -b(b+a)a -(2n+1)(b-a) - 2 b]\ ,
\]
\begin{equation}
\label{e.c.2}
J^0_n = \frac{1}{2} (b^2-a^2 - n-1)\ ,
\end{equation}
\[
J^-= \frac{b+a}{\sqrt{2}}\ ,
\]
where $n \in C$ . In this case the vacuum state
\begin{equation}
\label{e.c.3}
(b+a)|0>\  = \ 0 ,
\end{equation}
differs from (\ref{e1.2}).
If $n$ is a non-negative integer, then (\ref{e.c.2}) possesses
a finite-dimensional irreducible representation in the Fock space
\begin{equation}
\label{e.c.4}
{\cal P}_{n}(b) \ = \ \langle 1, (b-a), (b-a)^2, \dots , (b-a)^n \rangle \ ,
\end{equation}
of dimension $\dim{\cal P}_{n}=(n+1)$.

Taking $a,b$ in the realization (\ref{e1.3}) and substituting them into
(\ref{e.c.2}), we obtain
\[
J^+_n = \frac{1}{2^{3/2}}[x^3 + \pa_x^3 -x(x+\pa_x)\pa_x
-(2n+1)(x-\pa_x) - 2 \pa_x]\ ,
\]
\begin{equation}
\label{e.c.5}
J^0_n = \frac{1}{2} (x^2-\pa_x^2 - n-1)\ ,
\end{equation}
\[
J^-= \frac{x+\pa_x}{\sqrt{2}}\ ,
\]
which represents the $sl_2$-algebra by means of differential operators of finite
order (but not of first order as in (\ref{e.a.4})). The operator $J^0_n$ coincides
with the Hamiltonian of the harmonic oscillator (with the reference point
for eigenvalues changed). The vacuum state is
\begin{equation}
\label{e.c.6}
|0>\  = \ e^{-\frac{x^2}{2}}\ ,
\end{equation}
and the representation space is
\begin{equation}
\label{e.c.7}
{\cal P}_{n}(x) \ = \ \langle 1, x, x^2, \dots , x^n \rangle e^{-\frac{x^2}{2}}\ ,
\end{equation}
(cf.(\ref{e.c.4})).

{\bf (d).}
The following three operators
\[
J^+ = \frac{a^2}{2}\ ,
\]
\begin{equation}
\label{e.d.1}
J^0 = - \frac{\{ a,b\}}{4}\ ,
\end{equation}
\[
J^-= \frac{b^2}{2}\ ,
\]
are generators of the $sl_2$-algebra and the quadratic Casimir operator for
this representation is
\[
C_2 = \frac{3}{16}\ .
\]
This is the so-called metaplectic representation of $sl_2$ (see, for example,
\cite{per}). This representation is infinite-dimensional. Taking the
realization (\ref{e1.2}) or (\ref{e1.4}) of the Heisenberg algebra we get
the well-known representation
\begin{equation}
\label{e.d.2}
J^+ = \frac{1}{2}\pa_x^2\ ,\ J^0 = - \frac{1}{2}(x\pa_x-\frac{1}{2})\ ,
\ J^-= \frac{1}{2}x^2\
\end{equation}
in terms of differential operators, or
\[
J^+ = \frac{1}{2}{\cal D}_+^2\ ,\ J^0 = - \frac{1}{2}(x{\cal D}_--\frac{1}{2})\ ,
\]
\begin{equation}
\label{e.d.3}
J^-= \frac{1}{2}x(x-\de)(1-2\de{\cal D}_--\de^2{\cal D}_-^2)\ ,
\end{equation}
in terms of finite-difference operators, correspondingly.

{\bf (e).}
Take two operators $a$ and $b$ from the Clifford algebra $s_2$
\begin{equation}
\label{e.e.1}
             \{a,b\} \equiv ab  +  ba \ =\ 0\ ,\ a^2=b^2=1\ .
\end{equation}
Then the operators
\begin{equation}
\label{e.e.2}
J^1= a\ ,\ J^2= b\ , \ J^3= ab\ ,
\end{equation}
form the $sl_2$-algebra.

{\bf (f).}
Take the 5-dimensional Heisenberg algebra
\begin{equation}
\label{e.f.1}
  [a_i,b_j] \ =\ \de_{ij},\ i,j=1,2,\ldots,p \ ,
\end{equation}
where $\de_{ij}$ is the Kronecker symbol and $p=2$.
The operators
\begin{equation}
\label{e.f.2}
J^1= b_1a_2\ ,\ J^2= b_2a_1\ , \ J^3= b_1a_1 - b_2a_2\ ,
\end{equation}
form the $sl_2$-algebra. This representation is reducible. If (\ref{e.f.1}) is
given the coordinate-momentum representation
\begin{equation}
\label{e.f.3}
a_i\ =\ \frac{d}{dx_i} \equiv \pa_i\ ,\ b_i\ =\ x_i\ ,
\end{equation}
where $x \in {\bf C^2}$, the representation (\ref{e.f.2}) becomes the
well-known vector-field representation. The vacuum is a constant.
Finite-dimensional representations appear if a linear space
of homogeneous polynomials of fixed degree is taken.

\begin{center}
{\bf 2.} $sl_3$-algebra.
\end{center}
\renewcommand{\theequation}{A.2.{\arabic{equation}}}
\setcounter{equation}{0}

{\bf (a).}
Take the Fock space associated with the 5-dimensional
Heisenberg algebra (\ref{e.f.1}) with vacuum
\begin{equation}
\label{e2.a.1}
a_i|0>\  = \ 0\ ,\ i=1,2
\end{equation}
One can show that the following operators are the generators of the $sl_3$-algebra
\[
J^+_1\  =\ b_1(b_1 a_1 \  + \ b_2 a_2  - \ n)\ , \
J^+_2\  =\ b_2(b_1 a_1 \  + \ b_2 a_2  - \ n) \ ,
\]
\[
J^-_1\  =\ a_1\ ,\ J^-_2\  =\ a_2\ ,
\]
\[
J^0_{21}\ =\ b_2 a_1 \ ,\ J^0_{12}\  =\ b_1 a_2 \ ,
\]
\begin{equation}
\label{e2.a.2}
 J^0_1\  =\ b_1 a_1 - b_2 a_2\ , \
 J^0_2\  =\ b_1 a_1 + b_2 a_2\ - \ \frac{2}{3} n \ ,
\end{equation}
where $n$ is a complex number. If $n$ is a non-negative integer,
(\ref{e2.a.2}) possesses a finite-dimensional representation and its reprentation
space is given by the inhomogeneous polynomials of the degree not higher than $n$
in the Fock space:
\begin{equation}
\label{e2.a.3}
{\cal P}_n = \langle  b_1^{n_1} b_2^{n_2}\ |\
0 \leq (n_1+n_2) \leq n \rangle \ .
\end{equation}

In the coordinate-momentum representation (\ref{e.f.3}) the representation
(\ref{e2.a.2}) becomes
\[
J^+_1\  =\ x_1(x_1 \pa_1 \  + \ x_2 \pa_2  - \ n)\ , \
J^+_2\  =\ x_2(x_1 \pa_1 \  + \ x_2 \pa_2  - \ n) \ ,
\]
\[
J^-_1\  =\ \pa_1\ ,\ J^-_2\  =\ \pa_2\ ,
\]
\[
J^0_{21}\ =\ x_2 \pa_1 \ ,\ J^0_{12}\  =\ x_1 \pa_2 \ ,
\]
\begin{equation}
\label{e2.a.4}
 J^0_1\  =\ x_1 \pa_1 - x_2 \pa_2\ , \
 J^0_2\  =\ x_1 \pa_1 + x_2 \pa_2\ - \ \frac{2}{3} n \ ,
\end{equation}
where the vacuum $|0>=1$ and for non-negative integer $n$ the space of the
finite-dimensional representation is given by
\begin{equation}
\label{e2.a.5}
{\cal P}_n = \langle  x_1^{n_1} x_2^{n_2}\ |\
0 \leq (n_1+n_2) \leq n \rangle \ .
\end{equation}

{\bf (b).}
An important example of a quantum canonical transformation of
the 5-dimensional Heisenberg algebra (\ref{e.f.1}) is a generalization of
(\ref{e.b.1}) and has the form
\[
\hat a_i=\frac{(e^{\de_i a_i} -1)}{\de_i}\ ,
\]
\begin{equation}
\label{e2.b.1}
\hat b_i= b_i e^{ - \de_i a_i}\ ,\ i=1,2\ ,
\end{equation}
where $\de_{1,2}$ are complex numbers. Under this transformation the vacuum
remains the same (\ref{e2.a.1}). Finally, we are led to the following
representation of the $sl_3$-algebra
\[
J^+_1\  =\ \hat b_1(\hat b_1 \hat a_1 \  + \ \hat b_2 \hat a_2  - \ n)\ , \
J^+_2\  =\ \hat b_2(\hat b_1 \hat a_1 \  + \ \hat b_2 \hat a_2  - \ n) \ ,
\]
\[
J^-_1\  =\ \hat a_1\ ,\ J^-_2\  =\ \hat a_2\ ,
\]
\[
J^0_{21}\ =\ \hat b_2 \hat a_1 \ ,\ J^0_{12}\  =\ \hat b_1 \hat a_2 \ ,
\]
\begin{equation}
\label{e2.b.2}
 J^0_1\  =\ \hat b_1 \hat a_1 - \hat b_2 \hat a_2\ , \
 J^0_2\  =\ \hat b_1 \hat a_1 + \hat b_2 \hat a_2\ - \ \frac{2}{3} n \ ,
\end{equation}
As in previous case (a), for a non-negative integer $n$ the representation
(\ref{e2.b.2}) becomes finite-dimensional with the corresponding representation
space given by (\ref{e2.a.3}). We should mention that in the coordinate-momentum
representation the operators $\hat a,\hat b$ can be rewritten in terms of
finite-difference operators (\ref{e.b.3}-\ref{e.b.4}), ${\cal D}^{(x,y)}_{\pm}$
and, finally, the generators become
\[
J^+_1\ =\ x(1-\de_1{\cal D}^{(x)}_-)(x{\cal D}^{(x)}_-\  +\ y{\cal D}^{(y)}_-
- \ n)\ ,
\]
\[
J^+_2\ =\ y(1-\de_2{\cal D}^{(y)}_-)(x {\cal D}^{(x)}_-\ +\ y{\cal D}^{(y)}_-
- \ n)\ ,
\]
\[
J^-_1\ =\ {\cal D}^{(x)}_+\ ,\
J^-_2\ =\ {\cal D}^{(y)}_+\ ,
\]
\[
  J^0_{21}\ =\ y(1-\de_2{\cal D}^{(y)}_-){\cal D}^{(x)}_+\ ,
\ J^0_{12}\ =\ x(1-\de_1{\cal D}^{(x)}_-){\cal D}^{(y)}_+\ ,
\]
\begin{equation}
\label{e2.b.3}
 J^0_1\  =\ x {\cal D}^{(x)}_- \ - \ y {\cal D}^{(y)}_-\ , \
 J^0_2\  =\ x {\cal D}^{(x)}_- \ + \ y {\cal D}^{(y)}_-\ - \ \frac{2n}{3}\ .
\end{equation}

{\bf (c).}
Another representation of the $sl_3$-algebra is related to the 7-dimensional
Heisenberg algebra (\ref{e.f.1}) for $p=3$. The generators are
\[
J^+_1\  =\ -(b_2 - b_1 b_3)a_1-b_2 b_3 a_2-b_3^2 a_3+nb_3\ ,
\]
\[
J^+_2=-b_1(b_2 - b_1 b_3)a_1-b_2^2 a_2-b_2 b_3 a_3-m b_1 b_3+(n+m)b_2,
\]
\[
J^-_1\  =\ a_2\ ,\ J^-_2\  =\ a_3\ ,
\]
\[
J^0_{32}\ =\ a_1\ +\ b_3 a_2 \ ,
\ J^0_{23}\ =\ -b_1^2 a_1\ +\ b_2 a_3\ +\ m b_1,
\]
\[
 J^0_1\  =\ -b_1 a_1\ +\ b_2 a_2\ +\ 2 b_3 a_3\ -\ n\ ,
\]
\begin{equation}
\label{e2.c.1}
 J^0_2\  =\ 2b_1 a_1\ +\ b_2 a_2\ -\ b_3 a_3\ -\ m\ ,
\end{equation}
where $m, n$ are real numbers. In the coordinate-momentum representation
of Heisenberg algebra
the algebra (\ref{e2.c.1}) becomes the $sl_3$-algebra of first order
differential operators in the regular representation (on the flag manifold)
\[
J^+_1\  =\ -(y - x z)\pa_x\ -\ yz\pa_y\ -\ z^2\pa_z\ +\ nz\ ,
\]
\[
J^+_2=-x(y - xz)\pa_x-y^2\pa_y-yz\pa_z-m xz+(n+m)y\ ,
\]
\[
J^-_1\  =\ \pa_y\ ,\ J^-_2\  =\ \pa_z\ ,
\]
\[
J^0_{32}\ =\ \pa_x\ +\ z\pa_y\ ,
\ J^0_{23}\ =\ -x^2\pa_x\ +\ y\pa_z\ +\ m x,
\]
\[
 J^0_1\  =\ -x\pa_x\ +\ y\pa_y\ +\ 2z\pa_z\ -\ n\ ,
\]
\begin{equation}
\label{e2.c.2}
 J^0_2\  =\ 2x\pa_x\ +\ y\pa_y\ -\ z\pa_z\ -\ m\ .
\end{equation}

Using the realization (\ref{e2.b.1}) of the generators of the Heisenberg
algebra $H_7$ and the coordinate-momentum representation, a
realization of the $sl_3$-algebra emerges in terms of finite-difference
operators acting on $C^3$ functions, which is similar to (\ref{e2.b.3}).

\begin{center}
{\bf 3.} $gl_2 \ltimes {\bf C}^{r+1}$-algebra
\end{center}

\renewcommand{\theequation}{A.3.{\arabic{equation}}}
\setcounter{equation}{0}

Among the subalgebras of the (extended) universal enveloping algebra of the
Heisenberg algebra $H_5$ there is the 1-parameter family of non-semi-simple
algebras $gl_2 \ltimes {\bf C}^{r+1}$:
\[ J^1\  =\  a_1 \ , \]
\[ J^2\  =\ b_1 a_1\ -\ {n \over 3} \ ,\ J^3\  =\ b_2 a_2\
-\ {n\over {3r}} \ , \]
\[ J^4\  =\ b_1^2 a_1 \  +\ r b_1 b_2 a_2 \ - \ n b_1 \ ,\]
\begin{equation}
\label{e3.1}
 J^{5+k}\  = \ b_1^k a_2\ ,\ k=0,1,\dots, r\ ,
\end{equation}
where $r=1,2,\ldots$ and $n$ is a complex number.
Here the generators $J^{5+k},\ k=0,1,\dots, r$ span the $(r+1)$-dimensional
abelian subalgebra ${\bf C}^{r+1}$. If $n$ is a non-negative integer, the
finite-dimensional representation in the corresponding Fock space occurs,
\begin{equation}
\label{e3.2}
{\cal P}_n = \langle  b_1^{n_1} b_2^{n_2}\ |\
0 \leq (n_1+r n_2) \leq n \rangle \ .
\end{equation}

Taking the concrete realization of the Heisenberg algebra in terms of
differential or finite-difference operators in two variables similar to
(\ref{e1.2}) or (\ref{e1.4}) respectively, in the generators (\ref{e3.1})
we arrive at the $gl_2 \ltimes {\bf C}^{r+1}$-algebra realized as the algebra
of first-order differential operators
\footnote{This algebra acting on functions of two complex variables realized
by vector fields was found by Sophus Lie and, recently, it has been extended
to the algebra of first order differential operators \cite{gko}.}
or finite-difference operators, respectively.
\pagebreak
\begin{center}
{\bf 4.} $gl_k$-algebra.
\end{center}
\renewcommand{\theequation}{A.4.{\arabic{equation}}}
\setcounter{equation}{0}

The minimal Fock space where the $gl_k$-algebra acts is associated with
the $(2k-1)$-dimensional Heisenberg algebra $H_{2k-1}$. The explicit formulas
for the generators are given by
\[
J_i^-\ =\ a_i  \ ,   \quad i=2,3,\ldots , k \ ,
\]
\[
J_{i,j}^0\ =\ b_i J_j^-\ =\ b_i a_j \ , \quad i,j=2,3,\ldots , k \ ,
\]
\[
J^0 = n - \sum_{p=2}^k b_p a_{p}\ ,
\]
\begin{eqnarray}
\label{e4.1}
J_i^+ = b_i J^0\ , \quad i=2,3,\ldots, k \ ,
\end{eqnarray}
where the parameter $n$ is a complex number. The generators $J_{i,j}^0$
span the algebra $gl_{k-1}$. If $n$ is a
non-negative integer, the representation (\ref{e4.1}) becomes the
finite-dimensional representation acting on the space of polynomials
\begin{equation}
\label{e4.2}
V_n(t)\ =\ \mbox{span} \{ b_2^{n_2} b_3^{n_3} b_4^{n_4} \ldots
b_{k}^{n_{k}}\ |\ 0 \leq \sum n_i \leq n\}\ .
\end{equation}

Substituting the $a,b$-generators of the Heisenberg algebra in the
coordinate-momentum representation into (\ref{e4.1}) and  using
the vacuum, $|0>=1$, we get a representation of the $gl_k$-algebra in
terms of first-order differential operators (see, for example, \cite{btw})
\[
J_i^- = { \pa \over \pa x_i}  \ ,   \quad i=2,3,\ldots , k \ ,
\]
\[
J_{i,j}^0 = x_i J_j^-=x_i { \pa \over \pa x_j} \ , \quad
i,j=2,3,\ldots , k \ ,
\]
\[
J^0 = n - \sum_{p=2}^k x_p \frac{\pa}{\pa x_{p}} \ ,
\]
\begin{eqnarray}
\label{e4.3}
J_i^+ = x_i J^0\ , \quad i=2,3,\ldots, k \ ,
\end{eqnarray}
which acts on functions of $x \in \bf C^{k-1}$. One of the generators,
namely $J^{0} + \sum_{p=2}^{k} J_{p,p}^{0}$ is proportional to a constant
and, if it is taken out, we end up with the $sl_k$-subalgebra of the original
algebra. The generators $J_{i,j}^0$ form the $sl_{k-1}$-algebra of the vector
fields. If $n$ is a non-negative integer, the representation (\ref{e4.3})
becomes the finite-dimensional representation acting on the space of polynomials
\begin{equation}
\label{e4.4}
V_n(x)\ =\ \mbox{span} \{ x_2^{n_2} x_3^{n_3} x_4^{n_4} \ldots
x_{k}^{n_{k}}\ |\ 0 \leq \sum n_i \leq n\}\ .
\end{equation}
This representation corresponds to a Young tableau of one row with $n$ blocks
and is irreducible.

If the $a,b$-generators of the Heisenberg algebra are taken in the form of
finite-difference operators (\ref{e1.4}) and are inserted into (\ref{e4.1}),
the $gl_k$-algebra appears as the algebra of the finite-difference
operators:
\[
J_i^-\ =\ {\cal D}^{(i)}_+ \ ,   \quad i=2,3,\ldots , k \ ,
\]
\[
J_{i,j}^0=x_i(1-\de_i{\cal D}^{(i)}_-) J_j^-=x_i(1-\de_i{\cal D}^{(i)}_-)
{\cal D}^{(j)}_+,\ i,j=2,3,\ldots k,
\]
\[
J^0 = n - \sum_{p=2}^k x_p{\cal D}^{(p)}_-\ ,
\]
\begin{eqnarray}
\label{e4.5}
J_i^+ = x_i(1-\de_i{\cal D}^{(i)}_-) J^0\ , \quad i=2,3,\ldots, k \ ,
\end{eqnarray}
where ${\cal D}^{(i)}_{\pm}$ denote the finite-difference operators
(cf.(\ref{e1.4})) acting in the direction $x_i$.

\begin{center}
{\Large Lie Super-Algebras}
\end{center}
\renewcommand{\theequation}{S.{\arabic{equation}}}
\setcounter{equation}{0}

In order to work with superalgebras we must introduce the super Heisenberg
algebra. This is the $(2k+2r+1)$-dimensional algebra which contains the
$H_{2k+1}$-Heisenberg algebra (\ref{e.f.1}) as a subalgebra and also
the Clifford algebra $s_r$ :
\[
\{a^{(f)}_i, a^{(f)}_j\}\ =\ \{b^{(f)}_i, b^{(f)}_j\}\ =\ 0\ ,
\]
\begin{equation}
\label{s.1}
\{a^{(f)}_i, b^{(f)}_j\}\ =\ \de_{ij}\ ,\ i,j=1,2,\ldots,r\ ,
\end{equation}
as another subalgebra. There are two widely used realizations of the
Clifford algebra (\ref{s.1}):
\begin{itemize}
\item[(i)] The fermionic analogue of the coordinate-momentum representation
(\ref{e.f.3}):
\begin{equation}
\label{s.2}
a^{(f)}_i=\tha_i^{+}\ ,\ b^{(f)}_i\ =\ \tha_i\ ,\ i=1,2,\ldots,r\ ,
\end{equation}
or, differently,
\begin{equation}
\label{s.3}
a^{(f)}_i=\frac{\pa}{\pa \tha_i}\ ,\ b^{(f)}_i\ =\ \tha_i\ ,\
i=1,2,\ldots,r \ ,
\end{equation}
and
\item[(ii)] The matrix representation
\[
a^{(f)}_i\ =\
\underbrace{\si^0 \otimes \ldots \otimes \si^0}_{i-1} \otimes\ \si^+ \otimes
 \underbrace{{\bf 1} \otimes \ldots \otimes {\bf 1}}_{r-i}\ ,
\]
\begin{equation}
\label{s.4}
\ b^{(f)}_i\ =\ \underbrace{\si^0 \otimes \ldots \otimes \si^0}_{i-1} \otimes
\ \si^- \otimes \underbrace{{\bf 1} \otimes \ldots \otimes {\bf 1}}_{r-i}\ ,\
\ i=1,2,\ldots,r\ ,
\end{equation}
where the $\si^{\pm, 0}$ are Pauli matrices in standard notation,
\[
\si^+ =
\ \left( \begin{array}{cc}
0 & 1 \\
0 & 0
\end{array}  \right)
\ ,\ \si^-=
\ \left( \begin{array}{cc}
0 & 0 \\
1 & 0
\end{array}  \right)
\ ,\ \si^0 =
\ \left( \begin{array}{cc}
1 & 0 \\
0 & -1
\end{array}  \right) \ .
\]
\end{itemize}

In what follows we will consider the Fock space and also the realizations
of the superalgebras assuming that the Clifford algebra generators
are taken in the fermionic representation (\ref{s.3}) or the matrix
representation (\ref{s.4}).

\begin{center}

{\bf 1.} $osp(2,2)$-algebra.
\end{center}
\renewcommand{\theequation}{S.1.{\arabic{equation}}}
\setcounter{equation}{0}

Let us define a spinorial Fock space as a linear space of all 2-component
spinors with normal ordered polynomials in $a,b$ as components and with a
definition of the vacuum
\[
|0>\  =\
\ \left( \begin{array}{c}
|0>_1 \\
|0>_2
\end{array}  \right)
\]
such that any component is annihilated by the operator $a$:
\begin{equation}
\label{s1.1}
a|0>_i=0\ , \ i=1,2
\end{equation}

{\bf (a).}
Take the Heisenberg algebra (\ref{e1.1}). Then consider the following
two sets of $2 \times 2$ matrix operators:
\[
T^+ = b^2 a - n b + b \si^-\si^+,
\]
\begin{equation}
\label{s1.a.1}
 T^0 = ba - {n \over 2} +{1 \over 2} \si^-\si^+\  ,
\end{equation}
\[
T^- = a \ ,
\]
\[
J = -{ n \over 2} - {1 \over 2} \si^-\si^+ \ ,
\]
called bosonic (even) generators and
\begin{equation}
\label{s1.a.2}
Q \ = \ { \si^+ \brack b\si^+}\ , \
{\bar Q} \ = \ { (b a -n)\si^- \brack -a\si^-}\ ,
\end{equation}
called  fermionic (odd) generators.
The explicit matrix form of the even generators is given by:
\[
T^+ =
\left( \begin{array}{cc}
J^+_n & 0 \\
0 & J^+_{n-1}
\end{array}  \right),
T^0=
\left( \begin{array}{cc}
J^0_n & 0 \\
0 & J^0_{n-1}
\end{array}  \right),
T^- =
\left( \begin{array}{cc}
J^- & 0 \\
0 & J^-
\end{array}  \right),
\]
\[
 J=
\left( \begin{array}{cc}
-\frac{n}{2} & 0 \\
0 & -\frac{n+1}{2}
\end{array}  \right)\ ,
\]
and of the odd ones by
\[
Q_1 =
\left( \begin{array}{cc}
0 & 1 \\
0 & 0
\end{array}  \right)\ ,\
Q_2=
\left( \begin{array}{cc}
0 & b \\
0 & 0
\end{array}  \right)\ ,
\]
\begin{equation}
\label{s1.a.3}
\overline{Q}_1 =
\left( \begin{array}{cc}
0 & 0 \\
ba - n & 0
\end{array}  \right)\ ,\
 \overline{Q}_2=
\left( \begin{array}{cc}
0 & 0 \\
-a & 0
\end{array}  \right)\ ,
\end{equation}
where the $J^{\pm,0}_n$ are the generators of $sl_2$ given by (\ref{e.b.1}).

The above generators span the superalgebra $osp(2,2)$ with the commutation
relations:
\[
[T^0 , T^{\pm}]= \pm T^{\pm} \quad , \quad [T^+ , T^-]= -2 T^0 \quad , \quad
 [J , T^{\al}]=0  \quad , \al = \pm,0 \]
\[   \{ Q_1,\overline{Q}_2 \}  = - T^- \quad , \quad \{ Q_2, \overline{Q}_1 \}
= T^+ \ , \]
\[ {1 \over 2} (\{ \overline{Q}_1, Q_1\} + \{ \overline{Q}_2, Q_2 \})
= J \ ,
\
 {1 \over 2} (\{ \overline{Q}_1, Q_1\} - \{ \overline{Q}_2, Q_2 \}) =  T^0
\ , \]
\[ \{ Q_1,Q_1\}=\{ Q_2,Q_2\}=\{ Q_1,Q_2\}=0\ , \]
\[ \{ \overline{Q}_1,\overline{Q}_1\}=
\{ \overline{Q}_2,\overline{Q}_2\}=\{ \overline{Q}_1,\overline{Q}_2\}=0 \ ,\]
\[  [Q_1 , T^+]=Q_2 \ , \ [Q_2 , T^+]=0 \ , \ [Q_1 ,
T^-]=0 \ , \ [Q_2 , T^-]=-Q_1 \ ,\]
\[  [\overline{Q}_1 , T^+]=0 \ , \ [\overline{Q}_2 , T^+] = -
\overline{Q}_1 \ , \
[\overline{Q}_1 ,T^-] = \overline{Q}_2 \ , \ [\overline{Q}_2 ,
T^-]=0 \ ,\]
\[    [Q_{1,2} , T^0]=\pm {1 \over 2} Q_{1,2} \quad , \quad
 [\overline{Q}_{1,2} , T^0]=\mp {1 \over 2} \overline{Q}_{1,2}\ ,
\]
\begin{equation}
\label{s1.a.4}
  [Q_{1,2}, J] = - {1 \over 2} Q_{1,2} \quad , \quad
 [ \overline{Q}_{1,2}, J] = {1 \over 2} \overline{Q}_{1,2} \ .
\end{equation}

If, in the expressions (\ref{s1.a.1})--(\ref{s1.a.2}), the parameter
$n$ is a non-negative integer , then (\ref{s1.a.1})--(\ref{s1.a.2})
possess a finite-dimensional representation in the spinorial Fock space
\begin{equation}
\label{s1.a.5}
{\cal P}_{n,n-1} \ = \ \left \langle  \begin{array}{c}
1,b, b^2,\dots,b^{n}  \\
1, b, b^2, \dots,b^{n-1} \end{array} \right \rangle =
\ \left \langle  \begin{array}{c}
{\cal P}_{n} \\ {\cal P}_{n-1}
\end{array} \right \rangle \ .
\end{equation}

If we take a representation (\ref{e1.3}) of the Heisenberg algebra, the
generators (\ref{s1.a.1}) become  $2 \times 2$ matrix differential
operators, where the bosonic generators are \cite{shtur}
\[
T^+ = x^2 \pa_x - n x + x \si^-\si^+,
\]
\begin{equation}
\label{s1.a.6}
 T^0 = x \pa_x - {n \over 2} +{1 \over 2} \si^-\si^+\  ,
\end{equation}
\[
T^- = \pa_x \ ,
\]
\[
J = -{ n \over 2} - {1 \over 2} \si^-\si^+ \ ,
\]
and the fermionic generators
\begin{equation}
\label{s1.a.7}
Q \ = \ { \si^+ \brack x\si^+}\ , \
{\bar Q} \ = \ { (x \pa_x -n)\si^- \brack -\pa_x \si^-}\ .
\end{equation}
The finite-dimensional representation space for non-negative
integer values of the parameter $n$ in (\ref{s1.a.6}-\ref{s1.a.7})
becomes a linear space of 2-component spinors with polynomial components:
\begin{equation}
\label{s1.a.8}
{\cal P}_{n,n-1}(x) \ = \ \left \langle  \begin{array}{c}
1, x, x^2,\dots,x^{n}  \\
1, x, x^2, \dots,x^{n-1} \end{array} \right \rangle =
\ \left \langle  \begin{array}{c}
{\cal P}_{n}(x) \\ {\cal P}_{n-1}(x)
\end{array} \right \rangle
\end{equation}

{\bf (b).}
Taking the quantum canonical transformation (\ref{e.a.5}) and substituting
it into (\ref{s1.a.1}) we arrive at the $osp(2,2)$-algebra analogue of the
representation (\ref{e.b.2}) for the $sl_2$-algebra,
\[
T^+ =
(\frac{b}{\de}-1)b e^{ - \de a}(1 - n - e^{ - \de a}+ \si^-\si^+)\ ,
\]
\begin{equation}
\label{s1.b.1}
 T^0 =
\frac{b}{\de} (1 - e^{ - \de a}) -\frac{n}{2}+\frac{\si^-\si^+}{2}\ ,
\end{equation}
\[
T^- = {1 \over \de} ( e^{ \de a}-1) \ ,
\]
\[
J = -\frac{1}{2} - \frac{\si^-\si^+}{2} \ ,
\]
 and,
\begin{equation}
\label{s1.b.2}
Q \ = \ { \si^+ \brack b e^{ - \de a}\si^+}\ , \
{\bar Q} \ = \ { \frac{b - be^{ - \de a} -n}{\de}\si^- \brack
\frac{1 - e^{ \de a}}{\de}\si^-}\ .
\end{equation}

Taking for the generators $a,b$ the coordinate-momentum realization (\ref{e1.3}),
we obtain a representation of the algebra $osp(2,2)$ in terms of
finite-difference operators
\[
T^+ =
(\frac{x}{\de}-1)x e^{ - \de \pa_x}(1 - n - e^{ - \de \pa_x}+
 \si^-\si^+)\ ,
\]
\begin{equation}
\label{s1.b.3}
 T^0 =
\frac{x}{\de} (1 - e^{ - \de \pa_x}) -\frac{n}{2}+\frac{\si^-\si^+}{2}\ ,
\end{equation}
\[
T^- = {1 \over \de} ( e^{ \de \pa_x}-1) \ ,
\]
\[
J = -\frac{1}{2} - \frac{\si^-\si^+}{2} \ ,
\]
 and,
\begin{equation}
\label{s1.b.4}
Q \ = \ { \si^+ \brack x e^{ - \de \pa_x}\si^+}\ , \
{\bar Q} \ = \ { \frac{x - xe^{ - \de \pa_x} -n}{\de}\si^- \brack
\frac{1- e^{ \de \pa_x}}{\de}\si^-}\ .
\end{equation}
Or, in terms of the operators ${\cal D}_{\pm}$, their explicit matrix forms are
the following
\[
T^+ =
\left( \begin{array}{cc}
J^+_n & 0 \\
0 & J^+_{n-1}
\end{array}  \right),
T^0 =
\left( \begin{array}{cc}
x {\cal D}_- - \frac{n}{2} & 0 \\
0 & x {\cal D}_- - \frac{n-1}{2}
\end{array}  \right),
T^- =
\left( \begin{array}{cc}
{\cal D}_+ & 0 \\
0 & {\cal D}_+
\end{array}  \right),
\]
\[
 J=
\left( \begin{array}{cc}
-\frac{n}{2} & 0 \\
0 & -\frac{n+1}{2}
\end{array}  \right)\ ,
\]
for the bosonic generators and
\[
Q_1 =
\left( \begin{array}{cc}
0 & 1 \\
0 & 0
\end{array}  \right)\ ,\
Q_2=
\left( \begin{array}{cc}
0 & x(1-\de{\cal D}_-) \\
0 & 0
\end{array}  \right)\ ,
\]
\begin{equation}
\label{s1.b.5}
\overline{Q}_1 =
\left( \begin{array}{cc}
0 & 0 \\
x {\cal D}_- - n & 0
\end{array}  \right)\ ,\
 \overline{Q}_2=
\left( \begin{array}{cc}
0 & 0 \\
-{\cal D}_+ & 0
\end{array}  \right)\ ,
\end{equation}
for the fermionic generators, where the generator $J^+_n$ is given by
(\ref{e.b.6}).

{\bf (c).}
The super-metaplectic representation of the $osp(2,2)$-algebra can be easily
constructed and has the following form. The even generators are given by
\[
T^+ =
\left( \begin{array}{cc}
\frac{a^2}{2} & 0 \\
0 & \frac{a^2}{2}
\end{array}  \right),
T^0=
\left( \begin{array}{cc}
-\frac{\{a,b\}}{4} & 0 \\
0 & -\frac{\{a,b\}}{4}
\end{array}  \right),
T^- =
\left( \begin{array}{cc}
\frac{b^2}{2} & 0 \\
0 & \frac{b^2}{2}
\end{array}  \right),
\]
\[
 J\ =\
\left( \begin{array}{cc}
\frac{1}{4} & 0 \\
0 & -\frac{1}{4}
\end{array}  \right)\ ,
\]
while the odd ones are
\[
Q_1 =
\left( \begin{array}{cc}
0 & -\frac{b}{\sqrt{2}} \\
0 & 0
\end{array}  \right)\ ,\
Q_2=
\left( \begin{array}{cc}
0 & \frac{a}{\sqrt{2}} \\
0 & 0
\end{array}  \right)\ ,
\]
\begin{equation}
\label{s1.c.1}
\overline{Q}_1 =
\left( \begin{array}{cc}
0 & 0 \\
\frac{a}{\sqrt{2}} & 0
\end{array}  \right)\ ,\
 \overline{Q}_2=
\left( \begin{array}{cc}
0 & 0 \\
\frac{b}{\sqrt{2}} & 0
\end{array} \right)\ .
\end{equation}

Taking the realization of the Heisenberg algebra $H_3$ in terms of the
differential or finite-difference operators (\ref{e1.2}), (\ref{e1.4}),
respectively, and inserting it into (\ref{s1.c.1}) we end up with a
realization of the super-metaplectic representation of the $osp(2,2)$-algebra
in terms of differential or finite-difference operators.

\begin{center}
{\bf 2.} $gl(k+1,r+1)$-superalgebra.
\end{center}

\renewcommand{\theequation}{S.2.{\arabic{equation}}}
\setcounter{equation}{0}

One of the simplest representations of the $gl(k+1,r+1)$-superalgebra can
be written as follows
\[
T_i^- = a_i  \ ,   \quad i=1,2,\ldots,k \ ,
\]
\[
T_{i,j}^0=b_i T_j^-=b_i a_j\ ,\quad i,j=1,2,\ldots,k\ ,
\]
\[
T^0 = n - \sum_{p=1}^k b_p a_p-\sum_{p=1}^r {\tha}_p \frac{\pa}{\pa {\tha}_p}\ ,
\]
\begin{equation}
\label{s2.1}
T_i^+ = b_i T^0\ , \quad i=1,2,\ldots,k\ ,
\end{equation}
\[
\overline{Q}_i^- = \frac{\pa}{\pa {\tha}_i}\ ,\quad i=1,2,\ldots,r\ ,
\]
\[
\overline{Q}_i^+ = {\tha}_i T^0\ ,\quad i=1,2,\ldots,r\ ,
\]
\[
Q_{ij}^- = {\tha}_i T_j^-={\tha}_i a_j\ ,
\quad i,j=1,2,\ldots,r\ ,
\]
\[
Q_{ij}^+ = b_i \overline{Q}_j^-=b_i\frac{\pa}{\pa {\tha}_i}\ ,
\quad i=1,2,\ldots,k\ ,\ j=1,2,\ldots,r\ ,
\]
\[
J_{i,j}^0=\tha_i \overline{Q}_i^-=\tha_i \frac{\pa}{\pa \tha_j}\ ,
\quad i,j=1,2,\ldots,r\ ,
\]
These generators can be
represented by the following $(k+p) \times (k+p)$ matrix,
\begin{equation}
\label{s2.2}
\left( \begin{array}{ccc}
k \times k & | & p \times k \\
 BB        & | &    BF      \\
----        &   &    ----   \\
k \times p & | & p \times p \\
 FB        & | &     FF
\end{array} \right)\ ,
\end{equation}
where the notation $B(F)B(F)$ means the product of a bosonic operator $B$
(fermionic $F$) with a bosonic operator $B$ (fermionic $F$). Correspondingly, the
operators $T$ in (\ref{s2.1}) are of $BB$-type (mixed with $FF$-type), $J$ are
of $FF$-type, while the rest operators are of $BF$-type.
The algebra is defined by the (anti)commutation relations
\[
\{[ E_{IJ}, E_{KL}]\} =\de_{IL}E_{JK} \pm \de_{JK} E_{IL}\ ,
\]
where the generalized indices $I,J,K,L\ =\ B, F$. Anticommutators are taken for
generators of $FB, BF$ types only, while for all other cases the defining
relations are given by commutators. The dimension of the algebra is $(k+p)^2$.

The generators $J_{i,j}^0$ span the  $sl_k$-algebra of the vector fields.
The parameter $n$ in (\ref{s2.1}) can be any complex number. However, if $n$
is a non-negative integer, the representation (\ref{s2.1}) becomes the
finite-dimensional representation acting on a subspace of the Fock space
\begin{equation}
\label{s2.3}
V_n(b)=\mbox{span} \{ b_1^{n_1} b_2^{n_2} b_3^{n_3} \ldots b_k^{n_{k}}
\tha_1^{m_1}\tha_2^{m_2} \ldots \tha_r^{m_r}|0 \leq \sum n_i + \sum m_j \leq n\}.
\end{equation}

Taking the coordinate-momentum realization of the Heisenberg algebra
(\ref{e.f.3}) in the generators (\ref{s2.1}), we obtain the
$gl(k+1,r+1)$-superalgebra realized in terms of first order
differential operators (see, for example, \cite{btw}):
\[
T_i^- = \frac{\pa}{\pa x_i}  \ ,   \quad i=1,2,\ldots,k \ ,
\]
\[
T_{i,j}^0=x_i T_j^-=x_i \frac{\pa}{\pa x_j}\ ,\quad i,j=1,2,\ldots,k\ ,
\]
\[
T^0 = n - \sum_{p=1}^k x_p \frac{\pa}{\pa x_p}-\sum_{p=1}^r {\tha}_p
\frac{\pa}{\pa {\tha}_p}\ ,
\]
\begin{equation}
\label{s2.4}
T_i^+ = x_i T^0\ , \quad i=1,2,\ldots,k\ ,
\end{equation}
\[
\overline{Q}_i^- = \frac{\pa}{\pa {\tha}_i}\ ,\quad i=1,2,\ldots,r\ ,
\]
\[
\overline{Q}_i^+ = {\tha}_i T^0\ ,\quad i=1,2,\ldots,r\ ,
\]
\[
Q_{ij}^- = {\tha}_i T_j^-={\tha}_i \frac{\pa}{\pa x_j}\ ,
\quad i=1,2,\ldots,r;j=1,2,\ldots,k\ ,
\]
\[
Q_{ij}^+ = x_i \overline{Q}_j^-=x_i\frac{\pa}{\pa {\tha}_i}\ ,
\quad i=1,2,\ldots,k;j=1,2,\ldots,r\ ,
\]
\[
J_{i,j}^0=\tha_i \overline{Q}_i^-=\tha_i \frac{\pa}{\pa \tha_j}\ ,
\quad i,j=1,2,\ldots,r\ ,
\]
which acts on functions in ${\bf C}^k \otimes {\bf G}^r$.

A combination of the generators
$J^{0} + \sum_{p=1}^{k} T_{p,p}^{0}+\sum_{p=1}^{r} J_{p,p}^{0}$,
is proportional to a constant and, if it is taken out, we end up with the
superalgebra $sl(k+1,r+1)$. The generators
$T_{i,j}^0,\ J_{p,q}^0,\ i,j=1,2,\ldots,k\ ,\ p,q=1,2,\ldots,r$ span the algebra
of the vector fields $gl(k,r)$. The parameter $n$ in (\ref{s2.4}) can be any
complex number. If $n$ is a non-negative integer, the representation
(\ref{s2.1}) becomes the finite-dimensional representation acting on the
space of polynomials
\begin{equation}
\label{s2.5}
V_n(t)=\mbox{span} \{ x_1^{n_1} x_2^{n_2} x_3^{n_3} \ldots x_k^{n_{k}}
\tha_1^{m_1}\tha_2^{m_2} \ldots \tha_r^{m_r}|
0 \leq \sum n_i + \sum m_j \leq n\}.
\end{equation}
This representation corresponds to a Young tableau of one row with $n$ blocks
in the bosonic direction and is irreducible.

If the $a,b$-generators of the Heisenberg algebra are taken in the form of
finite-difference operators (\ref{e1.4}) and are inserted into (\ref{e4.1}),
the $gl(k+1,r+1)$-algebra appears as the algebra of the finite-difference
operators:
\[
T_i^- = {\cal D}^{(i)}_+  \ ,   \quad i=1,2,\ldots,k \ ,
\]
\[
T_{i,j}^0=x_i(1-\de_i{\cal D}^{(i)}_-) T_j^- =
x_i(1-\de_i{\cal D}^{(i)}_-) {\cal D}^{(j)}_+\ ,\quad i,j=1,2,\ldots,k\ ,
\]
\[
T^0 = n - \sum_{p=1}^k x_p{\cal D}^{(p)}_- -\sum_{p=1}^r {\tha}_p
\frac{\pa}{\pa {\tha}_p}\ ,
\]
\begin{equation}
\label{s2.6}
T_i^+ = x_i(1-\de_i{\cal D}^{(i)}_-) T^0\ , \quad i=1,2,\ldots,k\ ,
\end{equation}
\[
\overline{Q}_j^- = \frac{\pa}{\pa {\tha}_j}\ ,\quad j=1,2,\ldots,r\ ,
\]
\[
\overline{Q}_j^+ = {\tha}_j T^0\ ,\quad j=1,2,\ldots,r\ ,
\]
\[
Q_{ij}^- = {\tha}_i T_j^-={\tha}_i {\cal D}^{(j)}_+\ ,
\quad i=1,2,\ldots,r;j=1,2,\ldots,k\ ,
\]
\[
Q_{ij}^+ = x_i(1-\de_i{\cal D}^{(i)}_-) \overline{Q}_j^- =
x_i(1-\de_i{\cal D}^{(i)}_-)\frac{\pa}{\pa {\tha}_i},
\ i=1,2,\ldots k;j=1,2,\ldots r,
\]
\[
J_{i,j}^0=\tha_i \overline{Q}_i^-=\tha_i \frac{\pa}{\pa \tha_j}\ ,
\quad i,j=1,2,\ldots r\ ,
\]
It is worth mentioning that for the integer $n$, the algebra (\ref{s2.6})
has the same finite-dimensional representation (\ref{s2.5}) as the algebra
of the first order differential operators (\ref{s2.4}).

\begin{center}
{\Large Quantum Algebras}

$sl_{2q}$-algebra.
\end{center}
\renewcommand{\theequation}{Q.{\arabic{equation}}}
\setcounter{equation}{0}

Take two operators $\tilde a$ and $\tilde b$ obeying the commutation relation
\begin{equation}
\label{q.1.1}
             \tilde a \tilde b  -  q\tilde b\tilde a \ =\ 1\ ,
\end{equation}
with the identity operator on the r.h.s. They define the so-called
$q$-deformed Heisenberg algebra. Here $q \in C$. One can define a
$q$-deformed analogue of the universal enveloping algebra by taking all
ordered monomials ${\tilde b}^k{\tilde a}^m$. Introducing a vacuum
\begin{equation}
\label{q.1.2}
\tilde a|0>\  = \ 0\ ,
\end{equation}
in addition to the $q$-deformed analogue of the universal enveloping algebra
we arrive at a construction which is a $q$-analogue of Fock space.

It can be easily checked that the $q$-deformed Heisenberg algebra is a
subalgebra of the extended universal enveloping Heisenberg algebra.
This can be shown explicitly as follows. For any $q\in C$, two elements
of the extended universal enveloping Heisenberg algebra
\begin{equation}
\label{q.1.3}
\tilde a\ =\ \bigg(\frac{1}{b}\bigg)\bigg( {{q^{ba}-1}\over {q-1}}\bigg)\ ,
\ \tilde b\ =\ b\ ,
\end{equation}
obey the commutation relations (\ref{q.1.1}). It can be shown that the
universal enveloping Heisenberg algebra does not contain the $q$-deformed
Heisenberg algebra as a subalgebra. The formula (\ref{q.1.3}) allows us to
construct different realizations of the the $q$-deformed Heisenberg
algebra. One of them is a $q$-analogue of the coordinate-momentum
representation (\ref{e1.3}):
\begin{equation}
\label{q.1.4}
\tilde a\ =\ \tilde D_x\ ,\ \tilde b\ =\ x\ ,
\end{equation}
where
\begin{equation}
\label{q.1.5}
\tilde D_x f(x)\ =\ {{f(qx)-f(x)} \over {x(q-1)}}\ ,
\end{equation}
is the so-called Jackson symbol or the Jackson derivative.

Another realization of the (\ref{q.1.1}) appears if a quantum canonical
transformation of the Heisenberg algebra (\ref{e.b.1}) is taken:
\begin{equation}
\label{q.1.6}
\tilde a\ =\ \bigg(\frac{1}{b+\de}\bigg) e^{\de a}
\bigg( {{q^{\frac{b}{\de}(1-e^{-\de a})}-1}
\over {q-1}}\bigg)\ ,\ \tilde b\ =\ be^{ - \de a}\ ,
\end{equation}
where $\de$ is any complex number. In terms of translationally-covariant
finite-difference operators ${\cal D}_{\pm}$ the realization has the form
\begin{equation}
\label{q.1.7}
\tilde a\ =\ \bigg(\frac{1}{x+\de}\bigg)(\de{\cal D}_++1)\bigg( {{q^{x{\cal D}_-}-1}
\over {q-1}}\bigg)\ ,\ \tilde b\ =\ x(1-\de{\cal D}_-)\ .
\end{equation}
In these cases the vacuum is a constant, say, $|0>\  = \ 1$, as in the
non-deformed case.

The following three operators
\[
\tilde  J^+_{\al}\ =\ \tilde b ^2 \tilde a - \{ \al \} \tilde b \ ,
\]
\begin{equation}
\label{q.1.8}
\tilde  J^0_{\al}\ =\ \tilde  b\tilde  a - \hat{\al}\ ,
\end{equation}
\[
\tilde  J^-\  =\ \tilde a ,
\]
where $\{\al \} = {{1 - q^{\al}}\over {1 - q}}$  is so called $q$-number
and $\hat \al \equiv {\{\al \}\{\al+1\}\over \{2\al+2\}}$, are generators
of the $q$-deformed or quantum $sl_{2q}$-algebra. The operators
(\ref{q.1.8}) after multiplication by some factors, become
\[
\tilde  j^0 = {q^{-\al} \over q+1} {\{2\al+2\} \over \{\al+1\}}
\tilde J^0_\al \ ,
\]
\[
\tilde  j^{\pm} = q^{-\al/2} \tilde  J^{\pm}_{\al} \ ,
\]
and span the quantum algebra $sl_{2q}$ with the standard commutation relations
\cite{ot}\footnote{For discussion see \cite{at} as well},
\[
\tilde  j^0\tilde  j^+ \ - \ q\tilde  j^+\tilde  j^0 \ = \  \tilde  j^+\ ,
\]
\begin{equation}
\label{q.1.9}
 q^2 \tilde  j^+\tilde  j^- \ - \ \tilde  j^-\tilde  j^+ \
= \ - (q+1) \tilde  j^0 \ ,
\end{equation}
\[
q \tilde  j^0\tilde  j^- \ - \ \tilde  j^-\tilde  j^0 \
= \ - \tilde  j^- \ .
\]
\noindent
{\it Comment.}\  The algebra (\ref{q.1.9}) is known in literature as
{\it the second Witten quantum deformation} of $sl_2$ in the
classification of C.~Zachos \cite{z}).

In general, for the quantum $sl_{2q}$ algebra there are no polynomial Casimir
operators (see, for example, Zachos \cite{z}). However, in the representation
(\ref{q.1.8}) a relationship between generators analogous to the
quadratic Casimir operator appears
\[
q\tilde J^+_{\al}\tilde J^-_{\al} - \tilde J^0_{\al}
\tilde J^0_{\al} + (\{ {\al}+1 \}
- 2 \hat{\al}) \tilde J^0_{\al} = \hat{\al} (\hat{\al} - \{\al +1 \})\ .
\]
If $\al=n$ is a non-negative integer, then (\ref{q.1.8}) possesses a
finite-dimensional irreducible representation  in the Fock space
(cf.(\ref{e.a.2}))
\begin{equation}
\label{q.1.10}
{\cal P}_{n}(\tilde b) \ = \ \langle 1, \tilde b, \tilde b^2, \dots ,
\tilde b^n \rangle \ ,
\end{equation}
of the dimension $\dim{\cal P}_{n}=(n+1)$.

\end{document}